\documentclass[aos,MSNbibl,seceqn,dvips]{arximspdf}

\doi{10.1214/14-AOS1234} %kopijuoti is PTS
\volume{42}
\issue{6}
\pubyear{2014}
\firstpage{2141}
\lastpage{2163}
\docsubty{FLA}

\makeatletter
\newcommand{\rrvert}{\vert}
\newcommand{\rrVert}{\Vert}
\newcommand{\llvert}{\vert}
\newcommand{\llVert}{\Vert}
\newtheorem{lemma}{Lemma}[section]
\newtheorem{prop}[lemma]{Proposition}
\newtheorem{cor}[lemma]{Corollary}
\newtheorem{theorem}[lemma]{Theorem}
\newproclaim{remark}[lemma]{Remark}
\def\Chi{\raise.3ex \hbox{\large$\chi$}}
\def\cS{{\mathcal F}}
\def\cS{{\mathcal S}}
\def\diff{\,\triangle\,} %{\hbox{\tiny$\Delta$}}
\def\Prob{\operatorname{Prob}}
\def\exp{\operatorname{exp}}
\makeatother

\begin{document}
\begin{frontmatter}

\title{Classification algorithms using adaptive~partitioning\thanksref{T1}}
\runtitle{Classification algorithms using adaptive partitioning}

\begin{aug}
\author[A]{\fnms{Peter}~\snm{Binev}\ead[label=e1]{binev@math.sc.edu}},
\author[B]{\fnms{Albert}~\snm{Cohen}\corref{}\ead[label=e2]{cohen@ann.jussieu.fr}},
\author[C]{\fnms{Wolfgang}~\snm{Dahmen}\ead[label=e3]{dahmen@igpm.rwth-aachen.de}}
\and
\author[D]{\fnms{Ronald}~\snm{DeVore}\ead[label=e4]{rdevore@math.tamu.edu}}
\runauthor{Binev, Cohen, Dahmen and DeVore}
\affiliation{University of South Carolina,
Universit\'e Pierre et Marie Curie,\\
RWTH Aachen
and
Texas A\&M University}
\address[A]{P. Binev\\
Department of Mathematics\\
University of South Carolina\\
Columbia, South Carolina 29208\\
USA\\
\printead{e1}}
\address[B]{A. Cohen\\
Laboratoire Jacques-Louis Lions\\
UPMC Universit\'e Paris 06, UMR 7598\\
F-75005, Paris\\
France \\
\printead{e2}}
\address[C]{W. Dahmen\\
Institut f\"ur Geometrie\\
\quad und Praktische Mathematik\\
RWTH Aachen\\
Templergraben 55, D-52056, Aachen\\
Germany\\
\printead{e3}}
\address[D]{R. DeVore\\
Department of Mathematics\\
Texas A\&M University\\
College Station, Texas 77840\hspace*{29pt}\\
USA\\
\printead{e4}}
\end{aug}
\thankstext{T1}{Supported by the Office of Naval Research Contracts
ONR-N00014-08-1-1113, ONR-N00014-09-1-0107;
the AFOSR Contract FA95500910500; the ARO/DoD Contract W911NF-07-1-0185;
the NSF Grants DMS-09-15231, DMS-12-22390 and DMS-09-15104;
the Special Priority Program SPP 1324, funded by DFG;
the French--German PROCOPE contract 11418YB;
the Agence Nationale de la Recherche (ANR) project ECHANGE (ANR-08-EMER-006);
the excellence chair of the Foundation ``Sciences Math\'ematiques de
Paris'' held by
Ronald DeVore. This publication is based on work supported by Award No.
KUS-C1-016-04, made by King Abdullah University of Science and
Technology (KAUST).}

% HISTORY:
\received{\smonth{11} \syear{2013}}
\revised{\smonth{5} \syear{2014}}

% ABSTRACT
%
\begin{abstract}
Algorithms for binary classification based on adaptive tree
partitioning are formulated and
analyzed for both their risk performance and their friendliness to
numerical implementation.
The algorithms can be viewed as generating a set approximation to the
Bayes set and thus fall
into the general category of \textit{set estimators}.
In contrast with the most studied tree-based algorithms, which utilize
piecewise constant approximation
on the generated partition
[\textit{IEEE Trans. Inform. Theory} \textbf{52} (2006) 1335--1353;
\textit{Mach. Learn.} \textbf{66} (2007) 209--242], we consider decorated
trees, which allow us to derive higher order methods.
Convergence rates for these methods are derived in terms the parameter
$\alpha$ of margin conditions
and a rate $s$ of best approximation of the Bayes set by decorated
adaptive partitions.
They can also be expressed in terms of the Besov smoothness $\beta$ of
the regression function that governs
its approximability by piecewise polynomials on adaptive partition.
The execution of the algorithms does not require knowledge of the
smoothness or margin conditions.
Besov smoothness conditions are weaker than the commonly used H\"older
conditions, which govern approximation by nonadaptive partitions,
and therefore for a given regression function can result in a higher
rate of convergence.
This in turn mitigates the compatibility conflict between smoothness
and margin parameters.
\end{abstract}

% KEYWORDS
\end{frontmatter}

%s1 #&#
\section{Introduction}\label{sec1}

A large variety of methods has been developed for classification of
randomly drawn data.
Most of these fall into one of two basic categories: \textit{set
estimators} or \textit{plug-in estimators}.
Both of these families are based on some underlying form of approximation.
In the case of set estimators, one directly approximates the \textit{Bayes
set}, using elements from a
given family of sets. For plug-in estimators, one approximates the
underlying \textit{regression
function} and builds the classifier as a level set of this approximation.

The purpose of this paper is to introduce a family of classification
algorithms using tree-based adaptive partitioning and to analyze
their risk performance
as well as their friendliness to numerical implementation. These algorithms
fall into the category of set estimators.
Tree-based classification algorithms have been well studied since their
introduction
in \cite{CART}, and their convergence properties have been
discussed both in terms of oracle inequalities
and minimax convergence estimates; see, for example, \cite{SN} and
\cite{BSRM}. Among the specific features
of the approach followed in our paper are (i) the use of decorated
trees which allow us to derive faster rates for certain classes of
distributions than
obtained when using standard trees
and (ii) a convergence analysis based on nonlinear approximation theory
which allows us to significantly weaken the usual assumptions
that are made to establish a given convergence rate. More detailed
comparisons with existing methods and results are decribed later.

We place ourselves in the following setting of binary classification.
Let $X\subset\mathbb{R}^d$,
$Y=\{-1,1\}$ and $Z=X\times Y$. We assume that $\rho=\rho_X(x)\cdot
\rho(y|x)$ is a probability measure defined on $Z$. We denote by
$p(x)$ the probability that $y=1$ given $x$ and by $\eta(x)$ the
regression function
%
%
%e1.1 #&#
\begin{equation}
\label{reg} \eta(x):=\mathbb{E}(y|x)= p(x)-\bigl(1-p(x)\bigr)=2p(x)-1,
\end{equation}
where $\mathbb{E}$ denotes expectation. For any set $S$, we use the notation
%
%
%e1.2 #&#
\begin{equation}
\label{defs} \rho_S:=\rho_X(S)=\int
_S d\rho_X\quad\mbox{and}\quad\eta_S:=\int
_S \eta\,d\rho_X.
\end{equation}
A classifier returns
the value $y=1$ if $x$
%in
{is}
in some set $\Omega\subset X$ and $y=-1$ otherwise. Therefore, the
classifier is given by a function $T_\Omega=\Chi_\Omega-\Chi
_{\Omega^c}$ where $\Omega$ is some $\rho_X$ measurable set, and
$\Omega^c$ is its complement. With a slight abuse
of notation, we sometimes refer to the set $\Omega$ itself as the classifier.
We denote by $R(\Omega):= \mathbb{P}\{T_\Omega(x)\neq y\}$
the risk (probability of misclassification) of this classifier,
and by $\Omega^*:=\{x\dvtx \eta(x)\geq0\}$ the Bayes classifier
which minimizes this risk $R(\Omega)$, or equivalently, maximizes the
quantity $\eta_\Omega$ among all possible sets $\Omega$.

We measure the performance of a classifier $\Omega$
by the \textit{excess risk}
%
%
%e1.3 #&#
\begin{equation}
\label{erisk} R(\Omega)-R\bigl(\Omega^*\bigr)= \int_{\Omega\diff\Omega^*}
\llvert\eta\rrvert\,d\rho_X,
\end{equation}
with $A\diff B:=(A- B) \cup(B- A)$ the symmetric difference between
$A$ and $B$.

Given the data $\mathbf{z}=(z_i)_{i=1}^n$, $z_i=(x_i,y_i)$, $i=1,\ldots,n$,
drawn independently according to $\rho$,
a classification algorithm uses the draw to find a set $\widehat\Omega
=\widehat\Omega(\mathbf{z})$
to be used as a classifier. Obtaining a concrete estimate of the decay
of the excess risk for a given classification algorithm
as $n$ grows requires assumptions on the underlying measure $\rho$.
These conditions are usually spelled out by assuming that $\rho$ is in
a \emph{model class} ${\mathcal M}$.
Model classes are traditionally formed by two ingredients:
(i)~assumptions on the
behavior of $\rho$ near the boundary of the Bayes set $\Omega^*$
and~(ii)~assumptions on the smoothness of the regression function $\eta$.

Conditions that clarify (i) are called margin conditions and are
an item of many recent papers; see, for example, \cite{MN,Ts1}. One
formulation
(sometimes referred to as the Tsybakov condition) requires that
%
%
%e1.4 #&#
\begin{equation}
\label{tsy1} \rho_X\bigl\{x\in X\dvtx \bigl\llvert\eta(x)\bigr\rrvert
\le t\bigr\}\le C_\alpha t^\alpha, \qquad0<t\le1,
\end{equation}
for some constant $C_\alpha>0$ and $\alpha\geq0$. This condition
becomes more stringent as $\alpha$ tends to $+\infty$. The limiting case
$\alpha=\infty$, known as Massart condition,
means that for some $A>0$, we have $\llvert\eta\rrvert>A$ almost everywhere.
A common choice for (ii) in the classification literature (see, e.g.,
\cite{AT}) is that $\eta$ belongs to the H\"older space $C^\beta$.
This space can be defined for any $\beta>0$ as the set of functions
$f$ such that
%
%
%e1.5 #&#
\begin{equation}
\bigl\llVert\Delta_h^m f\bigr\rrVert_{L_\infty}
\leq C\llvert h\rrvert^\beta,\qquad h\in\mathbb{R}^d,
\end{equation}
where $\Delta_h^m$ is the $m$th power of the
finite difference operator defined by $\Delta_h f=f(\cdot+h)-f$, with
$m$ being
the smallest integer such that $m\geq\beta$. The H\"older class
may be viewed intuitively as the set of functions whose derivatives of
fractional order $\beta$ belong to $L_\infty$.

An important observation is that there is a conflict between margin and
smoothness
assumptions, in the sense that raising the smoothness $\beta$
limits the
range of $\alpha$ in the margin condition. For example, when $\rho_X$
is the Lebesgue measure
on $X$, it is easily checked that the constraint $\alpha\beta\leq1$
must hold
as soon as the Bayes boundary $\partial\Omega^*$ has nonzero
$(d-1)$-dimensional
Hausdorff measure.

An instance of a convergence result exploiting (i) and (ii),
Theorem~4.3 in \cite{AT}, is that under the assumption that
the density of $\rho_X$ with respect to the Lebesgue measure
is uniformly bounded, certain classifiers based on plug-in rules
achieve in expectation the rate
%
%
%e1.6 #&#
\begin{equation}
\mathbb{E}\bigl(R(\widehat\Omega)-R\bigl(\Omega^*\bigr)\bigr) \leq
Cn^{-((1+\alpha
)\beta)/((2+\alpha)\beta+d)}
\label{piholderrate}
\end{equation}
if the margin assumption holds with parameter $\alpha$, and if $\eta$
belongs to the H\"older class $C^\beta$.

The classification algorithms that we study in this paper have natural links
with the process of approximating the regression function using
piecewise constant or piecewise polynomials on adaptive partitions,
which is an instance of \emph{nonlinear approximation}.
It is well known in approximation theory that, when using nonlinear methods,
the
smoothness condition needed to attain a specified rate can be
dramatically weakened.
This state of affairs is reflected in the convergence results for our
algorithms that are
given in Theorems \ref{settheorem} and \ref{settheorem2}. These
results say
that with high probability (larger than $1-Cn^{-r}$ where $r>0$ can be
chosen arbitrarily
large), our classifiers achieve the rate
%
%
%e1.7 #&#
\begin{equation}
R(\widehat\Omega)-R\bigl(\Omega^*\bigr) \leq C \biggl(\frac{n} {\log n}
\biggr)^{-((1+\alpha)\beta)/((2+\alpha)\beta+d)} \label{treerate}
\end{equation}
if the margin assumption holds with parameter $\alpha$ and if $\eta$
belongs to the
Besov space $B^\beta_\infty(L_p)$ with $p>0$ such that $\beta p> d$.
This Besov space is defined
by the condition
%
%
%e1.8 #&#
\begin{equation}
\label{besov-def} \bigl\llVert\Delta_h^m f\bigr\rrVert
_{L_p} \leq C\llvert h\rrvert^\beta,\qquad h\in\mathbb{R}^d,
\end{equation}
with $m$ being an integer such that $m>\beta$
and may be viewed as the set of function whose derivatives of
fractional order $\beta$ belong to $L_p$.
Notice that the constraint $\beta p> d$ ensures that $B^\beta_\infty
(L_p)$ is compactly embedded in
$L_\infty$.
Therefore our algorithm achieves
the same rate as (\ref{piholderrate}), save for the logarithm, however,
with a significant
weakening on the smoothness condition imposed on the regression
function because of
the use of adaptive partitioning.
In particular, an individual regression
function may have significantly higher smoothness $\beta$ in this
scale of Besov spaces
than in the scale of H\"older spaces, resulting in a better rate
when using our classifier.

In addition, the weaker smoothness requirement
for a given convergence rate allows us to alleviate the conflict between
smoothness and margin conditions in the sense that the
constraint $\alpha\beta\leq1$ can be relaxed
% is no more needed
when using the space $B^\beta_\infty(L_p)$; see~(\ref{weaker}).
Let us also observe that our risk bound in (\ref{treerate}) holds in
the stronger
sense of high probability, rather than expectation, and that no
particular assumption (such as equivalence with Lebesgue measure) is
made on the density of $\rho_X$. Finally, let us stress
that our algorithms are numerically implementable and do not
require the a priori knowledge of the parameters $\alpha$ and $\beta$.

The distinction between Theorems \ref{settheorem} and \ref{settheorem2}
is the range of $\beta$ for which they apply.
Theorem \ref{settheorem} only applies to the range $\beta\le1$ and
can be seen as the
analog of using piecewise constant approximation on adaptive partition
for plug-in estimators.
On the other hand, Theorem \ref{settheorem2} applies for any $\beta
\le2$.
The gain in the range of $\beta$ results from the fact that the
algorithm uses decorated trees.
This corresponds to piecewise affine approximation for plug-in methods.
In principle, one can extend
the values of $\beta$ arbitrarily by using higher polynomial order
decorated trees.
However, the numerical implementation of such techniques becomes
more problematic and is therefore not considered in this paper. In the
regression
context, piecewise polynomial estimators on adaptive partitions have been
considered \mbox{in~\cite{A,BCDD2}}.

Set estimators aim at approximating the Bayes set $\Omega^*$
by elements $S$ from a family of sets $\cS$ in the sense
of the distance defined by the excess risk. Our approach to deriving the
risk bounds in Theorems \ref{settheorem} and \ref{settheorem2}
is by splitting this risk into
%
%
%e1.9 #&#
\begin{equation}
\label{errsplit} R(\widehat\Omega)-R\bigl(\Omega^*\bigr)=\bigl
(R(\widehat\Omega)-R(
\Omega_\cS)\bigr)+\bigl(R(\Omega_\cS)-R\bigl(\Omega^*
\bigr)\bigr),
\end{equation}
where $\Omega_\cS:=\operatorname{argmin} _{S\in\cS}
R(S)=\operatorname{argmax} _{S\in\cS}\eta_S$.
The two terms are positive, and are, respectively, referred
to as the estimation error and approximation error.
We bound in Section~\ref{sec2} the estimation error by
introducing a certain modulus, which is defined
utilizing any available uniform estimate between
the quantity $\eta_{S}-\eta_{\Omega_S}$ and its empirical
counterpart computed from the draw.
For set estimators based on empirical risk minimization, we show in
Section~\ref{sec3}
how margin conditions can be used to {bound} this modulus,
and therefore the {estimation error} term.

In Section~\ref{sec4}, we turn to estimates for the approximation term.
This analysis is based on the smoothness of $\eta$ and
the margin condition.
A typical setting when building set classifiers is
a nested sequence $(\cS_m)_{m\ge1}$ of families of sets, that is,
$\cS_m\subset\cS_{m+1}$,
where $m$ describes the complexity of $\cS_m$ in the sense of VC dimension.
The value of $m$ achieving between the optimal balance
between the estimation and approximation terms
depends on the parameters $\alpha$ and $\beta$,
which are unknown. A standard model selection procedure
is discussed in Section~\ref{sec5} that reaches this balance for a variety
of model classes ${\mathcal M}={\mathcal M}(\alpha,\beta)$ over a
range of $\alpha$ and $\beta$.

Many ingredients of our analysis of general classification
methods appear in earlier works; see, for example,
\cite{BBL,DGL}.
However, in our view, the organization of the material
in these sections helps clarify various issues concerning
the roles of approximation and {estimation
error bounds}.

In Section~\ref{sec6}, we turn to our proposed classification methods
based on
adaptive partitioning. We analyze their performance using the results from
the previous sections and obtain the aforementioned Theorems \ref{settheorem}
and \ref{settheorem2}.
The implementation and complexity of these algorithms
are discussed in Section~\ref{sec7}.

%%%%%%%%%%%%%%%%%%
%s2 #&#
\section{A general bound for the estimation error in set
estimators}\label{sec2}\label{Sgeneral}
%%%%%%%%%%%%%%%%%%%%%%%%%%%

In view of
$\Omega^*=\operatorname{argmax} _{\Omega
\subset X}\eta_\Omega$, if
$\hat\eta_S$ is any empirical estimator for $\eta_S$, a natural way
to select a classifier within $\cS$ is by
%
%
%e2.1 #&#
\begin{equation}
\label{empclassifier} \widehat\Omega:=\widehat\Omega_\cS:=\operatorname{argmax}
\limits
_{S\in
\cS}\hat\eta_S.
\end{equation}
One of the most common strategies for building $\hat\eta_S$ is by
introducing the empirical counterparts to (\ref{defs}),
%
%
%e2.2 #&#
\begin{equation}
\label{empiricaleta} \bar\rho_S:=\frac{1} n\sum
_{i=1}^n\Chi_S(x_i)\quad\mbox{and}\quad\bar\eta_S=\frac{1} n\sum
_{i=1}^ny_i\Chi_S(x_i).
\end{equation}
The choice $\hat\eta_S=\bar\eta_S$ is equivalent to minimizing
the empirical risk over the family~$\cS$, namely choosing
%
%
%e2.3 #&#
\begin{equation}
\widehat\Omega_\cS=\overline\Omega_\cS:=\operatorname
{argmin}\limits
_{S\in\cS} \overline R(S), \qquad\overline R(S):=\frac{1} n \#
\bigl\{ i\dvtx T_S(x_i)\neq y_i\bigr\}.
\end{equation}
However, other choices of $\hat\eta_S$
are conceivable, leading to other classifiers.

We give in this section a general method for bounding the {estimation error},
whenever we have an empirical estimator $\hat\eta_S$ for $\eta_S$,
with a bound of the form
%
%
%e2.4 #&#
\begin{equation}
\label{fundamental} \bigl\llvert\eta_S -\eta_{\Omega_\cS}-(\hat\eta
_S-\hat\eta_{\Omega_\cS
})\bigr\rrvert\le e_n(S),
\end{equation}
for each set $S\in\cS$. In the case where we use for $\hat\eta_S$
the set estimators $\bar\eta_S$ defined in (\ref{empiricaleta}),
we have the following bound.

%
%
%th2.1 #&#
\begin{theorem}
\label{VCtheorem}
For any sufficiently large constant $A>0$, the following holds. If $\cS
$ is a collection of $\rho_X$ measurable sets
$S\subset X$ with finite VC dimension $V:=V_\cS$, and if
%
%
%e2.5 #&#
\begin{equation}
\label{VCL11} e_n(S): = \sqrt{\rho_{S\diff\Omega_{\cS}} \varepsilon
_n}+\varepsilon_n,\qquad%\e_n:=\frac1 n\((r+1+AV)\log n)\).
{\varepsilon_n:=
A \max\{r+1,V\} \frac{\log n}n,}
\end{equation}
%
%For $A$ sufficiently large, we have the following:
where $r>0$ is arbitrary, then there is an absolute constant $C_0$ such
that for any $n\ge2$, with probability at least $1- C_0n^{-r}$ on the
draw $\mathbf{z}\in Z^n$, we have
%
%
%e2.6 #&#
\begin{equation}
\label{VCL22} \bigl\llvert\eta_S -\eta_{\Omega_\cS}-(\bar\eta_S-\bar
\eta_{\Omega_\cS
})\bigr\rrvert\le e_n(S),
\qquad S\in\cS.
\end{equation}
\end{theorem}

The techniques for proving this result are well known in
classification, but
we cannot find any reference that give the bounds in this theorem in
the above form,
and therefore we give its proof in the supplementary material~\cite{BCDD3}.

%The above theorem could be considered as a generalization of the
%similar result for finite collections $\cS$ of sets received through
%an application of %Bernstein's inequality. This result uses $\e_n:=
%probability $1-2n^{-r}$.
%}

%
%
%re2.2 #&#
\begin{remark}\label{finiteS}
The above theorem covers, in particular, the case where $\cS$ is a
\emph{finite} collection of sets,
since then trivially $V_\cS\leq\#\cS$. %this actually Theorem 1.3.6
%in \cite{DGL}.
Alternatively, in this case, a~straightforward\vspace*{1pt} argument using
Bernstein's inequality
yields the same result with the explicit expression
$\varepsilon_n:=\frac{10(\log(\#\cS)+r\log n)}{3n}$ and
probability at least $1-2n^{-r}$.
\end{remark}

To analyze the {{estimation error}} in classifiers, we
define the
following modulus:
%
%
%e2.7 #&#
\begin{equation}
\label{defmod} \omega(\rho,e_n):= \sup\biggl\{\int
_{S\diff\Omega_\cS}\llvert\eta\rrvert\dvtx S\in\cS\mbox{ and } \int
_{S\diff\Omega_\cS}\llvert\eta\rrvert\le3e_n(S) \biggr\}.
\end{equation}
Notice that the second argument $e_n$ is not a number but rather a set function.
In the next section, we discuss this modulus in some detail and bring
out its relation to other ideas used in classification, such as margin
conditions. For now, we use it to prove
the following theorem.\vspace*{-1pt}

%
%
%th2.3 #&#
\begin{theorem}
\label{th31}
Suppose that for each $ S\in\cS$, we have that (\ref{fundamental})
holds with probability $1-\delta$. Then with this same probability, we have
%
%
%e2.8 #&#
\begin{equation}
\label{want} R(\widehat\Omega_\cS)-R(\Omega_\cS) \leq\max
\bigl\{\omega(\rho,e_n), a\bigl(\Omega^*,\cS\bigr) \bigr
\}%,\qquad S\in\cS, %a_{
\end{equation}
with $a(\Omega^*,\cS):=R(\Omega_\cS)-R(\Omega^*)$ being the
approximation error from (\ref{errsplit}).\vspace*{-1pt}
%given by \iref{bias}.
\end{theorem}

\begin{pf}
We consider any data $\mathbf{z}$ such that (\ref{fundamental})
holds and prove that (\ref{want}) holds for such $\mathbf{z}$. Let
$S_0:=\Omega_{\cS}\setminus\widehat\Omega_\cS$ and
$S_1:= \widehat\Omega_\cS\setminus\Omega_\cS$ so that $S_0\cup
S_1=\widehat\Omega_\cS\diff\Omega_\cS$. Notice that, in contrast to
$\Omega_\cS$ and $\widehat\Omega_\cS$, the sets $S_0,S_1$ are
generally not in $\cS$.
We start from the equality
%
%
%e2.9 #&#
\begin{equation}
\label{start} R(\widehat\Omega_\cS)-R(\Omega_\cS)=
\eta_{\Omega_\cS}-\eta_{\widehat
\Omega_\cS} =\eta_{S_0}-
\eta_{S_1}.
\end{equation}
We\vspace*{1pt} can assume that $\eta_{S_0}-\eta_{S_1}> 0$, since otherwise we
have nothing to prove. From the definition of $\widehat\Omega_\cS$, we
know that
$ \hat\eta_{\Omega_\cS}- \hat\eta_{\widehat\Omega_\cS}\le0$.
Using this in conjunction with~(\ref{fundamental}), we obtain\vspace*{-1pt}
%
%
%e2.10 #&#
\begin{equation}
\label{implies} \eta_{S_0}-\eta_{S_1}= \eta_{\Omega_\cS}-
\eta_{\widehat\Omega_\cS} \le e_n(\widehat\Omega_\cS).
\end{equation}
In going further, we introduce the following notation. Given a set
$S\subset X$, we denote by $S^+:=S\cap\Omega^*$ and $S^-:=S\cap
(\Omega^*)^c$. Thus $\eta\ge0$ on $S^+$ and $\eta< 0$ on $S^-$.
Also $S=S^+\cup S^-$ and $S^+\cap S^-=\varnothing$.
Hence we can write
$\eta_{S_0}-\eta_{S_1}= A-B$, where $A:= \eta_{S_0^+}-\eta
_{S_1^-}\geq0$ and $B:=\eta_{S_1^+}-\eta_{S_0^-}\geq0$.
Note that\vspace*{1pt} $A,B\ge0$. We consider two cases.

\begin{longlist}[\textit{Case} 1.]
\item[\textit{Case} 1.] If $A\le2B$, then\vspace*{-1pt}
%
%
%e2.11 #&#
\begin{equation}
\label{start1} R(\widehat\Omega_\cS)-R(\Omega_\cS)= A-B
\leq B \le a\bigl(\Omega^*,{\cS}\bigr),
\end{equation}
where we have used the fact that $S_1^+\subset\Omega^*\setminus
\Omega_{\cS}$ and $S_0^-\subset\Omega_{\cS}\setminus\Omega^*$.
\end{longlist}

\begin{longlist}[\textit{Case} 2.]
\item[\textit{Case} 2.] If $A>2B$, then, by (\ref{implies}),\vspace*{-1pt}
%
%
%e2.12 #&#
\begin{equation}
\label{c2} \qquad\int_{\widehat\Omega_\cS\diff\Omega_\cS}\llvert\eta\rrvert
=A+B\le3A/2\le
3(A-B)=3(\eta_{S_0}-\eta_{S_1})\le3e_n(\widehat
\Omega_\cS).
\end{equation}
This means that $\widehat\Omega_\cS$ is one of the sets appearing in the
definition of $\omega(\rho,e_n)$, and~(\ref{want}) follows in this
case from the fact that\vspace*{-1pt}
\[
\eta_{S_0}-\eta_{S_1}=A-B\le\int_{\widehat\Omega_\cS\diff\Omega
_\cS}
\llvert\eta\rrvert\leq\omega(\rho,e_n).
\]
\end{longlist}\upqed%
\end{pf}

From Theorem \ref{th31}, we immediately obtain the following corollary.
%

%co2.4 #&#
\begin{cor}
\label{cor1}
Suppose that for each $ S\in\cS$, (\ref{fundamental})
\emph{holds} with probability $1-\delta$. Then with this same
probability we have
%
%
%e2.13 #&#
\begin{equation}
\label{want1} R(\widehat\Omega_\cS)-R\bigl(\Omega^*\bigr) \leq\omega(
\rho,e_n)+2 {a\bigl(\Omega^*,\cS\bigr)}. %a_{\cS}(\Omega^*).
\end{equation}
\end{cor}

%
%
%re2.5 #&#
\begin{remark}
\label{remlebesgue} The corollary does not impose any particular
assumption on $\rho$ and $\cS$,
apart from finite VC dimension.
For later comparisons with existing results, we briefly mention how
$e_n(S)$ can be sharpened under additional assumptions on~$\rho_X$.
Assume that $\rho_X$ is (or is equivalent to) the Lebesgue measure,
and for any arbitrary integer $l\geq1$, consider a uniform partition
${\mathcal Q}$ of $X=[0,1]^d$ into $m=l^d$ cubes of side length
$l^{-1}$, providing the collection $\cS$ of all sets $S$ that are
unions of cubes from ${\mathcal Q}$.
Then, defining
%
%
%e2.14 #&#
\begin{equation}
\label{newen} e_n(S):=\rho_{S\diff\Omega_{\cS}}\sqrt{
\varepsilon_n}\qquad\mbox{where } \varepsilon_n:=
\frac{8(r+1) m(1+\log n)}{3n},
\end{equation}
application of a union bound on Bernstein's inequality
to the random variable $y\Chi_{\Omega_\cS}(x)-y\Chi_S(x)$ gives
%
%
%e2.15 #&#
\begin{equation}
\label{pr21a} \mathbb{P}\bigl\{\bigl\llvert\eta_S-\bar\eta_S- (\eta
_{\Omega_\cS
}-\bar\eta_{\Omega_\cS}) \bigr
\rrvert\leq e_n(S)\dvtx S\in\cS\bigr\}\geq1- C n^{-r},
% \leq2 \exp\Big\{-\frac{3n\rho_S\o\PB{\e}_\WD{n}}{8}
% \Big\}= 2 \exp\Big\{-(r+1)\WD{d}{\WDn{m^d}}\rho_S(1+\log n).
% \Big\}.
\end{equation}
where $C$ is an absolute constant depending on $r$.
\end{remark}

%
%
%re2.6 #&#
\begin{remark}
\label{rempi}
Theorem \ref{th31} can be applied to any classification method
that is based on an estimation $\hat\eta_S$ of $\eta_S$,
once the bounds for $\llvert\eta_S -\eta_{\Omega_\cS}-(\hat\eta_S-\hat
\eta_{\Omega_\cS})\rrvert$
in terms of $e_n(S)$ have been established for all $S\in\cS$.
This determines $\omega(\rho,e_n)$ and thereby gives a bound for the
{estimation error}.
\end{remark}

%
%
%re2.7 #&#
\begin{remark}
\label{rem1}
The usual approach to obtaining bounds on the performance of
classifiers is to assume at the outset that the underlying measure
$\rho$ satisfies a margin condition. Our approach is motivated by the
desire to obtain bounds with no assumptions on $\rho$. This is
accomplished by introducing the modulus $\omega$. As we discuss in the
following section, a margin assumption allows one to obtain
an improved bound on $\omega$ and thereby recover existing results in
the literature.
Another point about our result is that we do not assume that the Bayes
classifier $\Omega^*$ lies in $\cS$. In some approaches, as discussed
in the survey \cite{BBL}, one first bounds the {estimation
error} under this assumption, and then later removes this assumption
with additional arguments
that employ margin conditions.
\end{remark}

%%%%%%%%%%%%%%
%s3 #&#
\section{Margin conditions}\label{sec3}\label{Smargin}
%%%%%%%%%%%%%%%

The modulus $\omega$ introduced in the previous section is not
transparent and, of course, depends on the set function $e_n(S)$.
However, as we now show, for the types of $e_n$ that naturally occur,
the modulus is intimately connected with margin conditions. Margin
assumptions are one of the
primary ingredients in obtaining estimates on the performance of
empirical classifiers.
The margin condition (\ref{tsy1}) recalled in the \hyperref[sec1]{Introduction} has the
following equivalent formulation: for any measurable set $S$,
we have
%
%
%e3.1 #&#
\begin{equation}
\rho_S\le C_\gamma\biggl( \int_S
\llvert\eta\rrvert\biggr)^\gamma, \qquad\gamma:= \frac{\alpha
}{1+\alpha}
\in[0,1] \label{tsy}
\end{equation}
for some constant $C_\gamma>0$ and $\gamma\in[0,1]$. This condition
is void when $\gamma=0$ and
becomes more stringent as $\gamma$ tends to $1$. The
case $\gamma=1$ gives Massart's condition.

In going further, we define ${\mathcal M}^\alpha$ as the set of all
measures $\rho$ such that $\rho_X$ satisfies~(\ref{tsy1}) or equivalently (\ref{tsy}), and we define
%
%
%e3.2 #&#
\begin{equation}
\label{Malpha} \llvert\rho\rrvert_{{\mathcal M}^\alpha}:=\sup_{0<t\le1}
t^{-\alpha} \rho_X\bigl\{ x\in X\dvtx \bigl\llvert\eta(x)\bigr
\rrvert\le t\bigr\}.
\end{equation}

We want to bring out the connection between the modulus $\omega$ and
the condition~(\ref{tsy}). In the definition of $\omega$ and its
application to bounds on the
{estimation error}, we assume that we have an empirical estimator for
which (\ref{fundamental}) holds with probability $1-\delta$. Notice
that this is only assumed to hold for sets $S\in\cS$ which is a
distinction with (\ref{tsy}).
We shall make our comparison first when $e_n $ is of the form
$e_n(S)=\sqrt{\varepsilon_n\rho_S}+\varepsilon_n$
as appears for set estimators in Theorem \ref{VCtheorem}.

We introduce the function
%
%
%e3.3 #&#
\begin{equation}
\label{phidef} \phi(\rho,t):=\sup_{ \int_S\llvert\eta\rrvert \le3(
t+\sqrt{t\rho_S})}\int
_S\llvert\eta\rrvert,\qquad0<t\le1,
\end{equation}
where now in this definition, we allow arbitrary measurable sets $S$
(not necessarily from $\cS$).
Under our assumption on the form of $e_n$, we have $\omega(\rho,e_n)\le
\phi(\rho,\varepsilon_n)$, and so the decay of $\phi$
gives us a bound on the decay of $\omega$.
We say that $\rho$ satsifies the $\phi$-condition of order $s>0$ if
%
%
%e3.4 #&#
\begin{equation}
\phi(\rho,t)\le C_0t^s, \qquad0<t\le1 \label{new}
\end{equation}
for some constants $C_0 $ and $s>0$.

%
%
%le3.1 #&#
\begin{lemma}
\label{nclemma}
Suppose\vspace*{1.5pt} that $\rho$ is a measure that satisfies (\ref{tsy1}) for a
given value
of $ 0\leq\alpha\leq\infty$. Then $\rho$ satisfies the $\phi
$-condition (\ref{new}) for $s=\frac{1+\alpha}{2+\alpha}$ with
$C_0$ depending only on $C_\alpha$ and $\alpha$. Conversely, if $\rho
$ satisfies the
$\phi$-condition\vspace*{2pt} with $s=\frac{1+\alpha}{2+\alpha}$ and a constant
$C_0>0$, then it satisfies (\ref{tsy1}) for the value
$\alpha$ with the constant $C_\alpha$ depending only on $s$ and $C_0$.
\end{lemma}

\begin{pf}
Suppose that $\rho$ satisfies (\ref{tsy1}) for $\alpha$
and constant $C_\alpha$, which equivalently means
that it satisfies (\ref{tsy}) for $\gamma:= \frac{\alpha}{1+\alpha
}$ and a constant\vspace*{2pt} $C_\gamma$.
To check that the $\phi$-condition is satisfied for $s=\frac{1+\alpha
}{2+\alpha}=\frac{1}{2-\gamma}$, we let $t\in(0,1]$ be fixed and
let $S$ be such that %PB $\int_S\left|\eta\right| \le3(\sqrt{t
$\int_S\llvert\eta\rrvert \le3(\sqrt{t\rho_S}+t)$. From (\ref{tsy}),
%
%
%e3.5 #&#
\begin{equation}
\label{know1} \rho_S\le C_\gamma\biggl(\int
_S\llvert\eta\rrvert\biggr)^\gamma\le
C_\gamma3^\gamma( \sqrt{t\rho_S}
+t)^\gamma.
\end{equation}
From this, one easily derives %
%
%
%e3.6 #&#
\begin{equation}
\label{derives} \rho_S \le Mt^{\gamma/(2-\gamma)},
\end{equation}
with a constant $M$ depending only on $C_\gamma$ and $\gamma$.
To see this, suppose to the contrary that for some {(arbitrarily
large)} constant $M$
%
%
%e3.7 #&#
\begin{equation}
\label{contrary} \rho_S > Mt^{\gamma/(2-\gamma)}.
\end{equation}
Rewriting (\ref{know1}) as
\[
\rho_S^{(2-\gamma)/(2\gamma)}\leq C_\gamma^{1/\gamma} 3 \bigl(
t^{1/2}+t\rho_S^{-1/2} \bigr),
\]
and using (\ref{contrary}) to estimate $\rho_S$ on both sides from
below, we obtain
\[
M^{(2-\gamma)/(2\gamma)}t^{1/2} \leq C_\gamma^{1/\gamma}3 \bigl(
t^{1/2}+M^{-1/2} t^{(4-3\gamma)/(4-2\gamma)} \bigr).
\]
Since $0<\gamma\le1$, we have $\frac{4-3\gamma}{4-2\gamma}\geq
\frac{1}2$, which yields
\[
t^{1/2}\leq M^{-(2-\gamma)/(2\gamma)}C_\gamma^{1/\gamma}3 \bigl(1+
M^{-1/2} \bigr) t^{1/2}.
\]
When $M$ is chosen large enough, we have $M^{-(2-\gamma)/(2\gamma
)}C_\gamma^{1/\gamma}3 (1+ M^{-1/2} )<1$
which is a contradiction thereby proving (\ref{derives}).

It follows from (\ref{know1}) and (\ref{derives}) that
%
%
%e3.8 #&#
\begin{equation}
\int_S\llvert\eta\rrvert\le3(t+\sqrt{t
\rho_S})\le3\bigl(t+ Mt^{1/(2-\gamma)}\bigr)\le C_0t^{1/(2-\gamma)},
\end{equation}
where $C_0$ depends on $C_\gamma$ and $\gamma$.
Taking now a supremum over all such sets $S$ gives
%
%
%e3.9 #&#
\begin{equation}
\phi(\rho,t)\le C_0t^s,\qquad s=\frac{1}{2-\gamma},
\end{equation}
which is the desired inequality.

We now prove the converse. Suppose that $\rho$ satisfies the $\phi
$-condition of order $s=\frac{1+\alpha}{2+\alpha}$ with constant
$C_0$. We want to show that
%
%
%e3.10 #&#
\begin{equation}
\label{show} \rho_X\bigl\{x\dvtx\bigl\llvert\eta(x)\bigr\rrvert\le
h\bigr\}\le C_\alpha h^{\alpha},\qquad0\le h\le1,
\end{equation}
with $ C_\alpha$ depending only on $s$ and $C_0$. As we noted before,
%% in \iref{tsy1},
this is equivalent to condition (\ref{tsy}) of order $\gamma= \frac
{\alpha}{\alpha+1}$.
To prove (\ref{show}), it is enough to prove
%
%
%e3.11 #&#
\begin{equation}
\label{show1} \rho_X\bigl\{x\dvtx h/2\le\bigl\llvert\eta(x)\bigr
\rrvert\le h\bigr\}\le{C_\alpha' } h^{\alpha
},
\qquad0<h\le1,
\end{equation}
since then (\ref{show}) follows easily by a summation argument.
We fix $h$ and define $S:=\{x\dvtx h/2\le\llvert\eta(x)\rrvert\le h\}
$ and $t: =
h^2\rho_{S}\in(0,1]$. Then we have
%
%
%e3.12 #&#
\begin{equation}
\int_S\llvert\eta\rrvert\le h\rho_{S}=
\sqrt{t\rho_{S}}.
\end{equation}
This means that $S$ is an admissible set in the definition of $\phi
(\rho,t)$ {in (\ref{phidef})}. Hence from the $\phi$-condition
{(\ref{new})}, we know
%
%
%e3.13 #&#
\begin{equation}
h\rho_{S}/2\le\int_S\llvert\eta\rrvert\le{\phi(\rho,t)} \le C_0 t^s=C_0
\bigl(h^2\rho_{S}\bigr)^{s}.
\end{equation}
In other words, we have
%
%
%e3.14 #&#
\begin{equation}
\rho_{S} \le( 2C_0)^{1/(1-s)}h^{(2s-1)/(1-s)}=
(2C_0)^{1/(1-s)}h^{\alpha},
\end{equation}
which completes the proof.
\end{pf}

%
%
%re3.2 #&#
\begin{remark}
The purpose of this section is to show the connection
of the modulus $\omega(\rho, e_n)$ with the existing and well-studied margin
conditions. However, the estimates for performance given in (\ref
{want1}) can be
applied without any specific assumption such as a margin condition,
which corresponds to the case $\gamma=0$.
One could also examine other types of bounds for $\phi(\rho,t)$ than
the power bound (\ref{new}) and obtain similar results.
\end{remark}

%%%%%%%%%%%%%%%%%%%%%%%%%%%%%%%%%
%s4 #&#
\section{Bounds for the approximation error \texorpdfstring{$a(\Omega^*,\mathcal
{F})$}{$a(Omega^*,\mathcal
{F})$}}\label{sec4}\label{Sbias}
%%%%%%%%%%%%%%%%%%%%

The approximation error $a(\Omega^*,\cS)$ depends on $\rho$ and the
richness of the collection $\cS$.
A typical setting starts with a nested sequence $(\cS_m)_{m\geq1}$ of
families of sets, that is such that $\cS_m\subset\cS_{m+1}$
for all $m\geq1$. The particular value of $m$ and the collection $\cS
_m$ that is used for a given draw of the data depends on $n$ and
properties of $\rho$ (such as the smoothness of $\eta$ and margin
conditions) and is usually chosen through some form of model selection
as discussed further. In order to analyze the performance of such
classification algorithms, we would like to know conditions on $\rho$
that govern the behavior of the approximation error as $m\to\infty$.
We study results of this type in this section.

The error
%
%
%e4.1 #&#
\begin{equation}
\label{errorm} a_m(\rho):=a\bigl(\Omega^*, \cS_m\bigr),
\qquad m\geq1,
\end{equation}
is monotonically decreasing. We define the approximation class
${\mathcal A}^s=\break  {\mathcal A}^s((\cS_m)_{m\geq1})$ as the set of all
$\rho$ for which
%
%
%e4.2 #&#
\begin{equation}
\label{approxclass} \llvert\rho\rrvert_{{\mathcal A}^s}:=\sup_{m\ge1}m^sa_m(
\rho)
\end{equation}
is finite.
Our goal is to understand what properties of $\rho$ guarantee
membership in~${\mathcal A}^s$. In this section, we give sufficient
conditions for $\rho$ to be in an approximation classes ${\mathcal
A}^s$ for both set estimators and plug-in estimators. These conditions
involve the smoothness (or approximability) of $\eta$ and margin conditions.

%that we have a monotone sequence $(\cS_m)_{m=1}^\infty$, where each $
Given a measure $\rho$, it determines the regression function $\eta$
and the Bayes set $\Omega^*:=\{x\dvtx\eta(x)>0\}$. We fix such a $\rho$,
and for each $t\in\mathbb{R}$, we define the level set $\Omega(t):=\{
x\dvtx \eta(x)\ge t\}$. Notice that $\Omega(t)\subset\Omega(t')$ if
$t\ge t'$. Also,
%
%
%e4.3 #&#
\begin{equation}
\label{notice} \bigl\{x\dvtx \bigl\llvert\eta(x)\bigr\rrvert< t\bigr
\}\subset
\Omega(-t)\setminus\Omega(t)\subset\bigl\{ x\dvtx \bigl\llvert\eta
(x)\bigr\rrvert\le
t\bigr\}.
\end{equation}

For each $m=1,2,\ldots,$ we define
%
%
%e4.4 #&#
\begin{equation}
\label{tm} t_m:=t_m(\rho,\cS_m):= \inf
\bigl\{t>0\dvtx \exists S\in\cS_m\mbox{ s.t. } \Omega(t)\subset
S\subset
\Omega(-t)\bigr\}.
\end{equation}
For convenience, we assume that there is always an $S_m^*\in\cS_m$
such that $\Omega(t_m)\subset S_m^*\subset\Omega(-t_m)$. [If no such
set exists then one replaces $t_m$ by $t_m+\varepsilon$ with
$\varepsilon> 0$ arbitrarily small and arrives at the same conclusion
(\ref{upperest}) given below.] It follows that
%
%
%e4.5 #&#
\begin{equation}
\label{contain} \Omega^*\diff S_m^*\subset\Omega(-t_m)
\setminus\Omega(t_m).
\end{equation}
If $\rho$ satisfies the margin condition (\ref{tsy1}), then
%
%
%e4.6 #&#
\begin{equation}
\label{upperest} a_m(\rho)\le\int_{\Omega^*\vartriangle S_m^*}\llvert
\eta\rrvert\,d\rho_X\le{C_\alpha} t_m\cdot
t_m^\alpha= {C_\alpha} t_m^{\alpha+1}.
\end{equation}
Thus a sufficient condition for $\rho$ to be in ${\mathcal A}^s$ is
that $t_m^{\alpha+1}\le Cm^{-s}$.

The following example illustrates how the margin condition (\ref{tsy1})
combined with
H\"older smoothness of the regression function
implies that $\rho$ belongs to the approximation class ${\mathcal
A}^s$ where $s$ depends on the margin and smoothness parameters.
To be specific, let
$X=[0,1]^d$. Let ${\mathcal D}$ be the collection of dyadic cubes $Q$
contained in $X$,
that is, cubes $Q\subset X$ of the form $Q=2^{-j}(k+[0,1]^d)$ with
$k\in\mathbb{Z}^d$ and $j\in\mathbb{Z}$.
Let ${\mathcal D}_j$, $j=0,1,\ldots,$ be the collection of dyadic cubes
of sidelength $2^{-j}$.
Let $\cS_{2^{dj}}$ be the collection of all sets of the form
$S_\Lambda=\bigcup_{Q\in\Lambda}Q$, where $\Lambda\subset{\mathcal
D}_j$. This corresponds to
the family $\cS$ considered in Remark \ref{remlebesgue} for
$m=2^{jd}$. In fact, $\#({\mathcal D}_j)=2^{jd}$ and $\#(\cS
_{2^{dj}})=2^{2^{jd}}$.
We complete the family $(\cS_m)_{m\geq1}$ by setting $\cS_m=\cS
_{2^{dj}}$ when $2^{dj}\leq m <2^{d(j+1)}$.
%

%
%
%pr4.1 #&#
\begin{prop}\label{proplebesgue}
We assume that $\rho$ has the two following properties:
\begin{longlist}[(ii)]
\item[(i)] the regression function $\eta$ is in the Lipschitz (or H\"older) space
$C^\beta$
for some $0<\beta\leq1$, that is,
\[
\llvert\eta\rrvert_{C^\beta}:=\sup\bigl\{\bigl\llvert\eta(x)-\eta(
\tilde{x})\bigr\rrvert\llvert x-\tilde{x}\rrvert^{-\beta}\dvtx
x,\tilde{x}
\in X\bigr\}<\infty;
\]
\item[(ii)]
$\rho$ satisfies the margin condition (\ref{tsy1}).

Then one has
%
%
%e4.7 #&#
\begin{equation}
\label{Asbeta} \rho\in{\mathcal A}^s={\mathcal A}^s
\bigl((\cS_{m})_{m\geq1}\bigr)\qquad\mbox{with } s:=
\frac{\beta(\alpha+1)} d.
\end{equation}
\end{longlist}
\end{prop}

\begin{pf}
We claim that
%
%
%e4.8 #&#
\begin{equation}
\label{ex1} a_{2^{{d}j}}(\rho)\le\bigl(M2^{-j\beta}
\bigr)^{\alpha+1},\qquad j\geq0,
\end{equation}
with $M:= 2^{-\beta} d^{\beta/2} \llvert\eta\rrvert_{C^\beta}$.
To this end, we first note that when $Q\in{\mathcal D}_j$, and $\xi
_Q$ is the center of $Q$, then %
%
%
%e4.9 #&#
\begin{equation}
\label{ex11} \bigl\llvert\eta(x)-\eta(\xi_Q)\bigr\rrvert\le M
2^{-j\beta}.
\end{equation}
We define $S_j\in\cS_{2^{dj}}$ as the union of all $Q\in{\mathcal
D}_j$ for which $\eta(\xi_Q)\ge0$. If $t:= M2^{-j\beta}$, then we
claim that
%
%
%e4.10 #&#
\begin{equation}
\label{contain1} \Omega(t)\subset S_{j}\subset\Omega(-t),\qquad j
\geq0.
\end{equation}
For example, if $x\in\Omega(t)$, then $\eta(x)\ge t$. So, if $x\in
Q$, then $\eta(\xi_Q)\ge0$ and hence $ Q\subset
S_j$. Similarly, if $x\in Q\subset S_j$, then $\eta(\xi_Q)\ge0$ and
hence $\eta(x)\ge-t$ for all $x\in Q$, and this implies the right
containment in (\ref{contain1}).
\end{pf}

It is well known that margin and smoothness conditions
are coupled, in the sense that higher values of $\alpha$ force the
regression function to have a sharper transition near the Bayes
boundary, therefore putting restrictions on its smoothness.
As an example, assume that $\rho_X$ is bounded from below
by the Lebesgue measure, that is, there exists a constant $c>0$ such
that for any $S\in\cS$,
\[
\rho_X(S) \geq c \llvert S\rrvert=c \int_S
dx.
\]
In the most typical
setting, the Bayes boundary $\partial\Omega^*$
is a ($d-1$)-dimensional surface of nonzero
${\mathcal H}^{d-1}$ Hausdorff measure.\vspace*{2pt} If $\eta\in C^\beta$
with $0\leq\beta\leq1$,
then $\llvert\eta(x)\rrvert$ is smaller than $t$ at any point
$x$ which is at distance less than $\llvert\eta\rrvert_{C^\beta
}^{1/\beta}
t^{1/\beta}$
from this boundary. It follows that
\[
\rho_X\bigl\{x\in X\dvtx \bigl\llvert\eta(x)\bigr\rrvert\leq t\bigr
\} \geq c_0t^{1/\beta},
\]
where\vspace*{1pt} $c_0$ depends on ${\mathcal H}^{d-1}(\partial\Omega^*)$
and $\llvert\eta\rrvert_{C^\beta}$, showing that $\alpha\beta\leq1$.
In such a case the approximation rate is therefore limited
by $s\leq\frac{1+\beta}{d}$.\vspace*{1pt}

As observed in \cite{AT} one can break this constraint
either by considering pathological examples,
such as regression functions
that satisfy ${\mathcal H}^{d-1}(\partial\Omega^*)=0$,
or by considering marginal measures $\rho_X$
that vanish in the vicinity of the Bayes boundary.
We show in Section \ref{Sadaptiveclass} that
this constraint can also be broken when
the H\"older spaces $C^\beta$ are %being
replaced by the Besov spaces $B^\beta_\infty(L_p)$, defined by (\ref
{besov-def}),
that govern the approximation rate when
$\cS_{2^{dj}}$ is replaced by a collection of adaptive partitions.

%%%%%%%%%%%%%%%%%%%%%%%%%%%%%
%s5 #&#
\section{Risk performance and model selection}\label{sec5}\label{Sriskperf}
%%%%%%%%%%%%%%%%%%%%%%%%%%%%%%%%%%%

In this section, we combine our previous {bounds for
approximation and
estimation errors} in order to obtain
an estimate for risk performance of classification schemes.

Let us assume that we have a sequence $(\cS_m)_{m\geq1}$ of families
$\cS_m$ of sets that are used to develop
a binary classification algorithm. We suppose that for some constant $C_0$,
%
%
%e5.1 #&#
\begin{equation}
\label{vcsm} VC(\cS_m)\leq C_0m,\qquad m\geq1,
\end{equation}
and we denote by $\overline\Omega_m$ the empirical risk minimization
classifier picked in $\cS_m$
according to (\ref{empclassifier}) with $\hat\eta_S=\bar\eta_S$.
{Theorem \ref{VCtheorem}} gives that such an estimator provides a
bound % \eref{VCL1}
{(\ref{VCL22})} with
\[
e_n(S)=\sqrt{ {\rho_{S\diff\Omega_{\cS}}} %\rho_S
\varepsilon_n
}+\varepsilon_n,\qquad% \e_n=Cm(\log n)/n
{\varepsilon_n=C
\frac{m\log n}{n}}
\]
and $C$ depending only on $r$ and $C_0$. If $\rho\in{\mathcal
A}^s((\cS_m)_{m\geq1})$, for some $s>0$, then according to Corollary
\ref{cor1}, for any $m\ge1$, we have with probability $1-n^{-r}$,
%
%
%e5.2 #&#
\begin{equation}
\label{risk11} R(\overline\Omega_m)-R\bigl( \Omega^*\bigr)\le\omega(
\rho,e_n)+2 \llvert\rho\rrvert_{{\mathcal A}^s}m^{-s}.
\end{equation}
If\vspace*{1pt} in addition $\rho$ satisfies the margin condition (\ref{tsy1}) of
order $\alpha>0$, then
using Lemma \ref{nclemma} and the fact that $\omega(\rho,e_n)\le
C\phi(\rho,\varepsilon_n)\le C\varepsilon_n^{(1+\alpha)/(2+\alpha)}$, we obtain
%
%
%e5.3 #&#
\begin{equation}
\label{risk2} R(\overline\Omega_m)-R\bigl( \Omega^*\bigr)\le C
\biggl(\frac{m\log
n}{n} \biggr)^{(1+\alpha)/(2+\alpha)}+2 \llvert\rho\rrvert
_{{\mathcal A}^s}m^{-s},
\end{equation}
where $C$ depends on $\llvert\rho\rrvert_{{\mathcal M}^\alpha}$.
If we balance the two terms appearing on the right in~(\ref{risk2}) by
taking $m=(\frac{n}{\log n})^{(1+\alpha)/((2+\alpha)s+1+\alpha
)}$, we obtain that, with probability at least $1-n^{-r}$,
%
%
%e5.4 #&#
\begin{equation}
R(\overline\Omega_m)-R\bigl( \Omega^*\bigr)\le C \biggl(
\frac{\log n}{n} \biggr)^{((1+\alpha)s)/((2+\alpha)s+1+\alpha)},
\label{risk3}
\end{equation}
where $C$ depends on $\llvert\rho\rrvert_{{\mathcal M}^\alpha}$ and
$\llvert\rho\rrvert_{A^s}$.
The best rates that one can obtain from the above estimate correspond to
$\alpha=\infty$ (Massart's condition) and $s\to\infty$ (the
regression function $\eta$ has arbitrarily high smoothness), and are
limited by the so-called fast rate ${\mathcal O}(\frac{\log n}{n})$.

To obtain the bound (\ref{risk3}), we need to know both $s$ and
$\alpha
$ in order to make the optimal choice of $m$ and $\cS_m$. Of course,
these values are not known to us, and to circumvent this we employ a
standard model selection technique based on
an independant validation sample.

\subsection*{Model selection}
Let $(\cS_m)_{m\ge1}$ be any
collection of set estimators.
For notational convenience, we assume that $n$ is even, that is,
$n=2\bar n$.
Given the draw $\mathbf{z}$, we divide $\mathbf{z}$ into two
independent sets
$\mathbf{z}'$ and $\mathbf{z}''$ of equal size $\bar n$.
\begin{longlist}[\textit{Step} 2.]
\item[\textit{Step} 1.] For each $1\le m\le\bar n$, we let $\overline
\Omega_m$ be defined by (\ref{empclassifier}) with $\cS=\cS_m$ and
$\mathbf{z}$ replaced by $\mathbf{z}'$.
\end{longlist}

\begin{longlist}[\textit{Step} 2.]
\item[\textit{Step} 2.] We\vspace*{2pt} now let $\overline\cS:=\{\overline\Omega_1,\ldots,\overline\Omega_{\bar n}\}$ and let $\widehat\Omega:=\overline
\Omega_{m^*}$ be the set chosen from~$\overline\cS$ by (\ref{empclassifier}) when using $\mathbf{z}''$.
\end{longlist}

The set $\widehat\Omega$ is our empirical approximation of $\Omega^*$
obtained by this model selection procedure. To see how well it
performs, let us now assume that $\rho\in{\mathcal A}^s$ and that
$\rho$ also satisfies the margin condition (\ref{tsy1}) for $\alpha$.
In {step~1}, we know that for each $m$, $ \overline\Omega_m$
satisfies (\ref{risk2}) with $n$ replaced by $\bar n$ with probability
at least $1-\bar n ^{-r}$. Thus, with probability $1-cn^{-r+1}$,
we have
%
%
%e5.5 #&#
%e5.6 #&#
\begin{eqnarray}\label{write1}
R(\overline\Omega_m)-R\bigl(\Omega^*\bigr) \le C
\biggl(m^{-s}+ \biggl(\frac
{m\log n}{n} \biggr)^{(1+\alpha)/(2+\alpha)} \biggr),
\nonumber\\[-8pt]\\[-8pt]
\eqntext{m=1,\ldots,\bar n.}
\end{eqnarray}
It follows that for $\overline\cS$ of {step~2}, we have
\[
a\bigl(\Omega^*,\overline\cS\bigr)=\min_{1\le m\le\bar n}\int
_{\overline\Omega_m\diff\Omega^*}\llvert\eta\rrvert\,d\rho_X \le C\min
_{1\le m\le\bar n} \biggl\{m^{-s}+ \biggl(\frac{m\log
n}{n}
\biggr)^{(1+\alpha)/(2+\alpha)} \biggr\}.
\]
Since $\#(\overline\cS)=\bar n=n/2$, we can take $\varepsilon
_n\le C\frac{\log n}{n}$ in Remark~\ref{finiteS} and a suitable
constant $C$ when
bounding performance on $\overline\cS$. Hence, from Corollary \ref
{cor1}, we have for the set $\widehat\Omega$ given by {step~2},
\begin{eqnarray*}
R(\widehat\Omega)-R\bigl(\Omega^*\bigr)& \le& 2a\bigl(\Omega^*,\overline
\cS\bigr) + C
\biggl(\frac{\log n}{n} \biggr)^{(1+\alpha)/(2+\alpha)}
\\
& \le& C\min_{1\le m\le\bar n} \biggl\{m^{-s}+ \biggl(
\frac
{m\log n}{n} \biggr)^{(1+\alpha)/(2+\alpha)} \biggr\}
\\
&&{}+ C \biggl
(\frac{\log n}{n}
\biggr)^{(1+\alpha)/(2+\alpha)}.
\end{eqnarray*}
In estimating the minimum, we choose $m$ that balances the two terms
and obtain
%
%
%e5.7 #&#
\begin{equation}
\label{risk31} R(\widehat\Omega)-R\bigl( \Omega^*\bigr)\le C \biggl(
\frac{\log n}{n} \biggr)^{((1+\alpha)s)/((2+\alpha)s+1+\alpha)}.
\end{equation}
Thus, the set $\widehat\Omega$, while not knowing $\alpha$ and $s$ gives
the same estimate we obtained earlier when assuming we knew $\alpha$
and $s$.

%
%
%re5.1 #&#
\begin{remark}
Note that we have done our model selection without
using a penalty term. The use of a penalty term would have forced us to
know the value of $\alpha$ in (\ref{tsy}). A discussion of why penalty
approaches may still be of interest in practice can be found in~\cite{BM}.
\end{remark}

A simple application of (\ref{risk31}) and Proposition \ref
{proplebesgue} gives an estimate in the general case.
In the case of Remark \ref{remlebesgue} one has $\omega(\rho,e_n)\le
C \varepsilon_n^{(\alpha+1)/2}$ and can balance the terms in
the estimate corresponding to (\ref{risk31}) by taking
$m:= (\frac{n}{(\log n)^{1/(2+d)}} )^{d/(2\beta+d)}$.
These give the following result.

%
%
%co5.2 #&#
\begin{cor}
\label{corprelim}
Let $(\cS_m)_{m\geq1}$
be the sequence of family of sets built from uniform partitions
as are used in Proposition \ref{proplebesgue}.
\begin{longlist}[(ii)]
\item[(i)] Assume that $\rho$ satisfies the margin condition (\ref
{tsy1}) of order $\alpha$ and that $\eta$ is H\"older continuous of
order $\beta$.
Then the classifier resulting from the above model selection satisfies
%
%
%e5.8 #&#
\begin{equation}
\label{betageneral} \Prob\biggl\{R(\overline\Omega_{m^*})-R\bigl(\Omega^*
\bigr)\le C \biggl(\frac
{\log n}{n} \biggr)^{((1+\alpha)\beta)/((2+\alpha)\beta
+d)} \biggr\}
\ge1-Cn^{-r},\hspace*{-25pt}
\end{equation}
where $C$ depends on $r,\llvert\rho\rrvert_{{\mathcal M}^\alpha},
\llvert\eta
\rrvert_{C^\beta}$.

\item[(ii)] If one assumes in addition that $\rho_X$ is equivalent to
the Lebesgue measure, one obtains
%
%
%e5.9 #&#
\begin{eqnarray}\label{betalebesgue}
&& \Prob\biggl\{R(\overline\Omega_{m^*})-R\bigl(\Omega^*
\bigr)\le C \biggl(\frac
{(\log n)^{1/(2+d)}}{n} \biggr)^{((1+\alpha)\beta)/((2+\alpha)\beta)} \biggr\}
\nonumber\\[-12pt]\\[-8pt]\nonumber
&&\qquad \ge1-Cn^{-r}.
\end{eqnarray}
\end{longlist}%
\end{cor}

Case~\textup{(ii)} illustrates the improvement of the rates that results from
constraining the marginal $\rho_X$. In view of our
earlier comments on the conflict between margin and H\"older smoothness
conditions of high order, the main deficiency
of both results is the underlying strong assumption of H\"older
smoothness. The sequence $(\cS_m)_{m\geq1}$
is based on uniform partitions and does not allow us to exploit weaker
Besov-smoothness conditions.
In what follows, we remedy this defect by turning to classification
algorithms based on \emph{adaptive} partitioning.
In doing so, we avoid any a priori constraints on $\rho_X$ and hence
use the set function $e_n$ given by (\ref{VCL11}).

%%%%%%%%%%%%%%%%%%%%%%%%%%%%%%%%%%%%%
%s6 #&#
\section{Classification using tree based adaptive partitioning}\label
{sec6}\label{Sadaptiveclass}
%%%%%%%%%%%%%%%%%%%%%%

One of the most natural ways to try to capture $\Omega^*$ is through
adaptive partitioning. Indeed, such partitioning
methods have the flexibility to give fine scale approximation near the
boundary of $\Omega^*$ but remain coarse away from the boundary. We
now give two examples. The first is based on simple dyadic tree
partitioning, while the second adds wedge ornation on the leaves of the
tree to enhance risk performance. For simplicity of presentation, we
only consider dyadic partitioning on the specific domain $X=[0,1]^d$,
even though our analysis covers far greater generality.

\subsection*{Algorithm \textup{I}: Dyadic tree partitioning}
We recall the dyadic cubes ${\mathcal D}$ introduced in Section \ref
{Sbias}. These cubes organize themselves into a tree with root $X$.
Each $Q\in{\mathcal D}_j$ has $2^d$ children which are its dyadic subcubes
from ${\mathcal D}_{j+1}$. A finite subtree ${\mathcal T}$ of
${\mathcal D}$ is a finite collection of cubes with the property that
the root $X$ is in ${\mathcal T}$, and whenever $Q\in{\mathcal T}$ its
parent is also in ${\mathcal T}$. We say a tree is \emph{complete} if,
whenever $Q$ is in ${\mathcal T}$, then all of its siblings are also in
${\mathcal T}$. The set ${\mathcal L}({\mathcal T})$ of leaves of such
a tree ${\mathcal T}$ consists of all the cubes $Q\in{\mathcal T}$
such that no child of $Q$ is in ${\mathcal T}$. The set of all such
leaves of a complete tree forms a partition of $X$.

Any finite complete tree is the result of a finite number of successive
cube refinements. We denote by $\frak T_m$ the collection of all
complete trees ${\mathcal T}$ that can be obtained using $m$
refinements. Any such tree ${\mathcal T}\in\frak T_m$ has $(2^d-1)m+1$ leaves.
We can bound the number of trees in ${\mathcal T}\in\frak T_m$ by assigning
a bitstream that encodes, that is, precisely determines, ${\mathcal T}$
as follows. Let ${\mathcal T}\in\frak T_m$. We order the children of~$X$ lexicographically and assign a one to every child which is refined
in ${\mathcal T}$ and a zero otherwise. We now consider the next
generation of cubes (i.e., the grandchildren of~$X$) in ${\mathcal T}$.
We know these grandchildren from the bits already assigned. We arrange
the grandchildren lexicographically and again assign them a one if they
are refined in ${\mathcal T}$ and a zero otherwise. We continue in this
way and receive a bitstream which exactly determines ${\mathcal T}$.
Since ${\mathcal T}$ has exactly $2^dm+1$ cubes, every such bitstream
has length $2^dm$ and has a one in exactly $m-1$ positions. Hence we have
%
%
%e6.1 #&#
\begin{equation}
\label{count1} \#(\frak T_m)\le \pmatrix{2^d m\cr m-1} \le
\frac{(2^d m)^{m}}{(m-1)!} \le e^{m}2^{dm}.
\end{equation}

For each ${\mathcal T}\in\frak T_m$ and any $\Lambda\subset{\mathcal
L}( {\mathcal T})$, we define $S=S_\Lambda:=\bigcup_{Q\in
\Lambda}Q$. We denote by $\cS_m$ the collection of all such
sets $S$ that can be obtained from a ${\mathcal T}\in\frak T_m$ and
some choice of $\Lambda$. Once ${\mathcal T}$ is chosen, there are
$2^{\#({\mathcal L}({\mathcal T}))}\le2^{2^dm}$ choices for $\Lambda
$. Hence
%
%
%e6.2 #&#
\begin{equation}
\label{count2} \#(\cS_m)\le a^m
\end{equation}
with $a:=e2^{d+2^d}$.

Given our draw $\mathbf{z}$, we use the set estimator and model selection
over $(\cS_m)_{m\ge1}$ as described in the previous section. We
discuss the numerical implementation of this algorithm in Section \ref
{SSnumerical}. This results in a set $\overline\Omega(\mathbf{z})$, and
we have the following theorem for its performance.
%

%th6.1 #&#
\begin{theorem}\label{settheorem}
\textup{(i)} For any $r>0$, there is a constant $c>0$ such that the
following holds. If $\rho\in{\mathcal A}^s$, $s>0$ and $\rho$
satisfy the margin condition {(\ref{tsy1})}, %\eref{tsy1},
then with probability greater than $ 1-cn^{-r+1}$, we have
%
%
%e6.3 #&#
\begin{equation}
\label{settheorem1} R\bigl(\overline\Omega(\mathbf{z})\bigr)-R\bigl
(\Omega^*\bigr)\le
C \biggl(\frac{\log
n}{n} \biggr)^{((1+\alpha)s)/((2+\alpha)s+1+\alpha)}
\end{equation}
with $C$ depending only on $d,r$, $\llvert\rho\rrvert_{{\mathcal
A}^s}$ and the
constant in (\ref{tsy1}).

\textup{(ii)} If $\eta\in B_\infty^\beta(L_p(X))$ with $0<\beta\le1$
and $p>d/\beta$ and if {$\rho$} %$\eta$
satisfies the margin condition {(\ref{tsy1})}, %\eref{tsy1},
then with probability greater than $ 1-cn^{-r+1}$, we have
%
%
%e6.4 #&#
\begin{equation}
\label{settheorem11} %PB{settheorem1}
R\bigl(\overline\Omega(\mathbf{z})\bigr)-R\bigl(
\Omega^*\bigr)\le C \biggl(\frac{\log
n}{n} \biggr)^{((1+\alpha)\beta)/((2+\alpha)\beta+d)},
\end{equation}
with $C$ depending only on $d,r$, $\llvert\eta\rrvert_{B_\infty^\beta
(L_p(X))}$
and the constant in (\ref{tsy1}).
\end{theorem}

\begin{pf}
Since $\log(\#(\cS_m))\le C_0m$ where $C_0$ depends only
on $d$, we have that $R(\Omega(\mathbf{z}))-R(\Omega^*)$ is bounded
by the
right-hand side of (\ref{risk31}) which proves (i).
We can derive (ii)\vspace*{1pt} from (i) if we prove that the assumptions on $\rho$
in (ii) imply that \mbox{$\rho\in{\mathcal A}^s$}, $s=\frac{(\alpha
+1)\beta}{d}$. To see that this is the case, we consider the
approximation of
$\eta$ by piecewise constants subordinate to partitions ${\mathcal
L}({\mathcal T})$, ${\mathcal T}\in\frak T_m$. It is known (see~\cite
{DNL}) that the Besov space assumption on $\eta$ implies that there is
a tree ${\mathcal T}_m$ and piecewise constant $\eta_m$ on ${\mathcal
L}({\mathcal T}_m)$ that satisfies $\llVert \eta-\eta_m\rrVert
_{L_\infty}\le
\delta_m= C_1\llvert\eta\rrvert_{B_\infty^\beta(L_p)}m^{-\beta/d}$
with $C_1$
depending on $p, \beta$ and $d$. Let $\Lambda:=\{Q\in{\mathcal
L}({\mathcal T}_m)\dvtx \eta_m(x)>0, x\in Q\}$ and
$\Omega_m:=\bigcup_{Q\in\Lambda_m}Q$. Then $\Omega
_m\in\cS_m$ and $\Omega_m\diff\Omega^*\subset\{x\dvtx \llvert\eta
(x)\rrvert\le
\delta_m\}$, and so
%
%
%e6.5 #&#
\begin{equation}
\label{settheorem21} %PB{settheorem2}
a_m(\rho)\le\int_{\Omega_m\diff\Omega^*}
\llvert\eta\rrvert\,d\rho_X\le C_\alpha\delta_m^{\alpha+1}
\le C_\alpha\bigl( C_1\llvert\eta\rrvert_{B_\infty^\beta(L_p)}
\bigr)^{\alpha+1}m^{-s},
\end{equation}
as desired.
\end{pf}

\subsection*{Algorithm \textup{II}: Higher order methods via decorated trees}
We want to remove the restriction $\beta\le1$ that appears in Theorem
\ref{settheorem} by enhancing the family of sets
$\cS_m$ of the previous section. This enhancement can be accomplished
by choosing, for each $Q\in{\mathcal L}({\mathcal T})$, a subcell of
$Q$ obtained by a hyperplane cut (henceforth called an \emph{H-cell})
and then taking a union of such subcells. To describe this, we note that,
given a dyadic cube $Q$, any ($d-1$)-dimensional hyperplane $H$
partitions $Q$ into at most two disjoint sets $Q^H_0$ and $Q^H_1$ which
are the intersections of $Q$ with the two open half spaces generated by
the hyperplane cut. By convention we include $Q\cap H$ in $Q^H_0$.
Given a tree ${\mathcal T}\in\frak T_m$, we denote by $\zeta
_{\mathcal T}$ any mapping that assigns to each $Q\in{\mathcal
L}({\mathcal T})$ an H-cell $\zeta_{\mathcal T}(Q)$. Given such a
collection $\{\zeta_{\mathcal T}(Q)\}_{Q\in{\mathcal L}({\mathcal T})}$,
we define
\[
S:=S({\mathcal T},\zeta):=\bigcup_{Q\in{\mathcal L}({\mathcal T})}
\zeta_{\mathcal T}(Q).
\]
For any given tree ${\mathcal T}$, we let $\cS_{\mathcal T}$ be the
collection of all such sets that result from arbitrary choices of
$\zeta$.
For any $m\ge1$, we define
%
%
%e6.6 #&#
\begin{equation}
\label{Sm} \cS_m:=\bigcup_{{\mathcal T}\in\frak T_m}
\cS_{\mathcal T}.
\end{equation}
Thus any such $S\in\cS_m$ is the union of H-cells of the $Q\in
{\mathcal L}(T)$, with one H-cell chosen for each $Q\in{\mathcal
L}({\mathcal T})$.
Clearly $\cS_m$ is infinite, however, the following lemma shows that
$\cS_m$ has finite VC dimension.

%
%
%le6.2 #&#
\begin{lemma}
\label{VClemma1}
If $\Gamma_1,\ldots,\Gamma_N$ are each collections of sets from $X$
with VC dimension $\le k$, then the collection $\Gamma:=\bigcup
_{i=1}^N\Gamma_i$ has VC dimension not greater than
$\max\{8\log N, 4k\}$.
\end{lemma}

\begin{pf}
We follow the notation of Section~9.4 in \cite{GKKW}. Let
us consider any set of points $p_1,\ldots, p_L$ from $X$. Then from
Theorem~9.2 in \cite{GKKW},
the shattering number of $\Gamma$ for this set of point satisfies
\[
s\bigl(\Gamma_j,\{p_1,\ldots, p_L\}\bigr)
\le\sum_{i=0}^k \pmatrix{L\cr i}=:\Phi(k,L)
\]
and therefore
\[
s\bigl(\Gamma,\{p_1,\ldots, p_L\}\bigr)\le N \Phi(k,L).
\]
By Hoeffding's inequality, if $k\leq L/2$, we have $2^{-L}\Phi
(k,L)\leq\exp(-2L\delta^2)$
with $\delta:=\frac{1} 2-\frac{k} L$. It follows that if $L>\max\{8\log
N, 4k\}$, we have
\[
s\bigl(\Gamma,\{p_1,\ldots, p_L\}\bigr)<2^L
N\exp(-L/8)<2^L,
\]
which shows that $\mathrm{VC}(\Gamma)\leq\max\{8\log N, 4k\}$.
\end{pf}

We\vspace*{1.5pt} apply Lemma \ref{VClemma1} with the role of the $\Gamma_j$ being
played by the collection $\cS_{\mathcal T}$, %PB $\cS_\cT:$,
${\mathcal T}\in\frak T_m$.
As shown in (\ref{count1}), % of the previous section,
we have $N= \#(\frak T_m)\leq e^m 2^{dm}$. We note next that the VC
dimension of
each $\cS_{\mathcal T}$ is given by
%
%
%e6.7 #&#
\begin{equation}
\label{VCST1} \mathrm{VC}(\cS_{\mathcal T})= (d+1)\#\bigl({\mathcal L}({
\mathcal T})\bigr)\le(d+1)2^d m.
\end{equation}
In\vspace*{1pt} fact, given ${\mathcal T}$ placing $d+1$ points in every $Q\in
{\mathcal L}({\mathcal T})$ shows that $(d+1)\#({\mathcal L}({\mathcal
T}))$ points can be shattered
since $d+1$ points can be shattered by hyperplanes in $\mathbb{R}^d$.
No matter how one distributes more than
$(d+1)\#({\mathcal L}({\mathcal T}))$ points in $X$,\vadjust{\goodbreak} at least one $Q\in
{\mathcal L}({\mathcal T})$ contains more than $d+1$ points. These points
can no longer be shattered by a hyperplane which confirms (\ref
{VCST1}). Lemma \ref{VClemma1} now says that
%
%
%e6.8 #&#
\begin{equation}
\label{VCSm} \mathrm{VC}(\cS_m)\leq\max\bigl\{ 8(d+2) m,
4(d+1)2^d m\bigr\} =C_d m,
\end{equation}
where $C_d:= \max\{ 8(d+2), 4(d+1)2^d\}$.

Given our draw $\mathbf{z}$, we use the set estimator and model
selection as
described in
Section~\ref{Sriskperf} %\ref{SSmodelsel}
with $\cS_m$ now given by (\ref{Sm}). This results in a set
$\overline\Omega(\mathbf{z})$, and we have the following theorem for
the performance of this estimator.
%

%th6.3 #&#
\begin{theorem}
\label{settheorem2}
\textup{(i)} For any $r>0$, there is a constant $c>0$ such that the
following holds. If $\rho\in{\mathcal A}^s$, $s>0$ and $\rho$
satisfy margin condition (\ref{tsy1}), then with probability greater
than $ 1-cn^{-r+1}$, we have
%
%
%e6.9 #&#
\begin{equation}
\label{settheorem12} R\bigl(\overline\Omega(\mathbf{z})\bigr)-R\bigl
(\Omega^*\bigr)\le
C \biggl(\frac{\log
n}{n} \biggr)^{((1+\alpha)s)/((2+\alpha)s+1+\alpha)}
\end{equation}
with $C$ depending only on $d,r$, $\llvert\rho\rrvert_{{\mathcal
A}^s}$ and the
constant in (\ref{tsy1}).

\textup{(ii)} If $\eta\in B_\infty^\beta(L_p(X))$ with $0<\beta\le2$
and $p> d/\beta$ and if {$\rho$} %$\eta$
satisfies the margin condition (\ref{tsy}), then with probability
greater than $ 1-cn^{-r+1}$, we have
%
%
%e6.10 #&#
\begin{equation}
\label{settheorem13} R\bigl(\overline\Omega(\mathbf{z})\bigr)-R\bigl
(\Omega^*\bigr)\le
C \biggl(\frac{\log
n}{n} \biggr)^{((1+\alpha)\beta)/( (2+\alpha)\beta+d)},
\end{equation}
with $C$ depending only on $d,r$, $\llvert\eta\rrvert_{B_\infty^\beta
(L_p(X))}$
and the constant in (\ref{tsy1}).
\end{theorem}

\begin{pf} In view of (\ref{VCSm}) we can invoke Theorem~\ref
{VCtheorem} %Lemma \ref{VClemma}
with
$\varepsilon_n = Cm\log n/\break n$, where $C$ depends on $d$ and $r$, to
conclude that $e_n(S) = \sqrt{\rho_{S\diff\Omega_{\cS
_m}}\varepsilon_n}+\varepsilon_n$ satisfies~(\ref{VCL22})
and hence
is an admissible set function for the modulus (\ref{defmod}).
Now (i) follows from
(\ref{risk31}).

To derive (ii) from (i), we prove that the assumptions on $\rho$ in
(ii) imply that $\rho\in{\mathcal A}^s$, $s=\frac{(\alpha+1)\beta}{d}$,
for $\beta\in(0,2]$.
To\vspace*{1pt} see that this is the case, we consider the approximation of
$\eta$ by piecewise \emph{linear} functions subordinate to partitions
${\mathcal L}({\mathcal T})$, ${\mathcal T}\in\frak T_m$. It is known
(see \cite{CDDD}) that the Besov space assumption on $\eta$ implies
that there is
a tree ${\mathcal T}_m$ and a piecewise {linear} function $\eta_m$ on
${\mathcal L}({\mathcal T}_m)$ that satisfies $\llVert \eta-\eta
_m\rrVert
_{L_\infty}\le\delta_m= C_1\llvert\eta\rrvert_{B_{\infty}^\beta
(L_p(X))}m^{-\beta/d}$.
Now for any cube $Q$ consider the H-cell mapping $\zeta_{\mathcal
T}(Q):= \{x\in Q\dvtx \eta_m(x)\geq0\}$. Then
\[
\Omega_m:=\bigcup_{Q\in{\mathcal L}({\mathcal T})}
\zeta_{\mathcal T}(Q)
\]
is in $\cS_m$ and $\Omega_m\diff\Omega^*\subset\{x\dvtx \llvert\eta
(x)\rrvert\le
\delta_m\}$ so that
%
%
%e6.11 #&#
\begin{equation}
\label{settheorem22} \qquad a_m(\rho)\le\int_{\Omega_m\diff\Omega^*}\llvert
\eta\rrvert\,d\rho_X\le C_\alpha\delta_m^{\alpha+1}
\le C_\alpha\bigl( C_1\llvert\eta\rrvert_{B_\infty^\beta(L_p)}
\bigr)^{\alpha+1}m^{-s},
\end{equation}
as desired.
\end{pf}

%
%
%re6.4 #&#
\begin{remark}
It is, in theory, possible to further extend
the range of $\beta$ by
considering more general decorated trees,
where for each considered cube $Q$,
we use an algebraic surface $A$ of degree
$k>1$ instead of a hyperplane $H$ that corresponds
to the case $k=1$. The
resulting families $\cS_m$ consist
of level sets of piecewise polynomials of degree $k$ on
adaptive partitions obtained by $m$ splits.
From this one easily
shows that the corresponding VC dimension is
again controlled by $m$ (with multiplicative constants
now depending both on $d$ and $k$) and that
(\ref{settheorem13}) now holds for all $ 0<\beta\leq k+1$.
However, the practical implementation
of such higher order classifiers appears to be difficult.
\end{remark}

We have seen in Section~\ref{sec5} that the approximation rate
for nonadaptive partitioning is also given by $s=\frac{\beta(\alpha+1)}{d}$,
but\vspace*{1pt} with $\beta$ denoting the smoothness of $\eta$
in the sense of the H\"older space $C^\beta$.
The results established in this section show that
the same approximation rate is obtained
under the weaker constraint that
$\eta\in B_\infty^\beta(L_p)$ with $p> d/\beta$
if we use adaptive partitioning.

We also observed in Section~\ref{sec5} that the
H\"older smoothness $\beta$ and the parameter $\alpha$
in the margin condition are coupled,
for example, by the restriction
$\alpha\beta\leq1$ when $\rho_X$ is bounded
from below by the Lebesgue measure.
Replacing the H\"older space $C^\beta$
by a Besov space $B_\infty^\beta(L_p)$ with $p> d/\beta$
allows us to relax the above constraint.
As a simple example consider the case
where $\rho_X$ is the Lebesgue measure and
\[
\eta(x)=\eta(x_1,\ldots,x_d)= \operatorname{sign}(x_1-1/2)
\llvert x_1-1/2\rrvert^\delta,
\]
for some $0<\delta\leq1$, so that
$\Omega^*=\{x\in X\dvtx x_1>1/2\}$,
and margin condition (\ref{tsy1}) holds with $\alpha$ up to $1/
\delta$.
Then one checks that $\eta\in B_\infty^\beta(L_p)$ for $\beta$ and
$p$ such
that $\beta\leq\delta+1/p$. The constraint $1/p< \beta/d$ may then be
rewritten as $\beta(1-1/d)<\delta$ or equivalently
%
%
%e6.12 #&#
\begin{equation}
\label{weaker} \alpha\beta(1-1/d)<1,
\end{equation}
which is an improvement over $\alpha\beta\leq1$.

%%%%%%%%%%%%%%%%%%%%%%%%%%%%%%%%

%%%%%%%%%%%%%%%%%%%%%%%%%%%%%%%%
%s7 #&#
\section{Numerical implementation}\label{sec7}\label{SSnumerical}
%%%%%%%%%%%%%%%%%%%%%%%%%%%%%%%

The results we have presented thus far on adaptive partitioning do not
constitute a numerical algorithm since we have not discussed how one
would find the sets $\overline\Omega_m\in\cS_m$ %given in
required by
(\ref{empclassifier}) and used in the model selection. We discuss this
issue next.

%%%%%%%%%%%%%%%%%%%%PB
%%%%%%%%%%%%%%%%%%%% green start

Given the draw $\mathbf{z}$, we consider the collection of all dyadic cubes
in ${\mathcal D}_0\cup\cdots\break \cup{\mathcal D}_{\bar n}$
with $\bar n=n/2$ which contain an $x_i$, $i=1,\ldots, \bar n$.
{These cubes form a tree ${\mathcal T}'(\mathbf{z})$ which we call the
\emph{occupancy tree}.}
Adding to all such cubes
their siblings, we obtain a complete tree ${\mathcal T}(\mathbf{z})$
whose leaves form a partition of $X$.
%which we call the {\it occupancy tree}. There are
%at most $n$ occupied cubes in $\cL(\cT(\bz))$
%and therefore $\#(\cL(\cT(\bz)))\leq2^d n$. It is easily checked that
%the cardinality of a tree does not exceed twice that of its leaves and
%therefore
%$\#(\cT(\bz))\leq2^{d+1}n$.

Let us first discuss the implementation of Algorithm I. For each %cube
%$Q\in\cT(\bz)$, and any
complete subtree ${\mathcal T}\subset{\mathcal T}(\mathbf{z})$ %
we define
%
%
%e7.1 #&#
\begin{equation}
\label{energy7.1} \gamma_{{\mathcal T}}:=\sum_{Q\in{\mathcal L}({\mathcal T})}
\max(\bar\eta_Q,0),
\end{equation}
which we call the \emph{energy} in ${\mathcal T}$. The set estimator
$\overline\Omega_m$ corresponds to a complete tree $\overline
{\mathcal T}_m\in\frak T_m$ which maximizes the above energy.
Note that several different trees may attain the maximum. Since only
the values $m=1,\ldots,\bar n$ are considered in the model
selection procedure, and since there is no gain
in subdividing a nonoccupied cube,
a maximizing tree is always a subtree of ${\mathcal T}(\mathbf{z})$.

Further, for each cube $Q\in{\mathcal T}(\mathbf{z})$, we denote by
$\frak
T_m(Q)$ the collection of all
complete trees ${\mathcal T}$ with root $Q$ obtained using at most $m$
subdivisions
and being contained in ${\mathcal T}(\mathbf{z})$. We then define
%
%
%e7.2 #&#
\begin{equation}
\label{energy1} \gamma_{Q,m}= \max_{{\mathcal T}\in\frak T_m(Q)}
\gamma_{\mathcal T}.
\end{equation}
Again, this maximum may be attained by several trees in $\frak T_m(Q)$.
In fact, if, for instance, for
a maximizer ${\mathcal T}\in\frak T_m(Q)$, $\bar\eta_R>0$
holds for all $R\in{\mathcal C}(R')\subset{\mathcal L}({\mathcal
T})$, the children of
some parent node $R'\in{\mathcal T}$, then the subtree $\widetilde
{{\mathcal T}}$ of ${\mathcal T}$ obtained by removing ${\mathcal
C}(R')$ from ${\mathcal T}$,
has the same energy. We denote by
${\mathcal T}(Q,m)$ any tree in $\frak T_m(Q)$ that attains the maximum
$\gamma_{Q,m}$. By convention, we set
%
%
%e7.3 #&#
\begin{equation}
{\mathcal T}(Q,m)=\varnothing,
\end{equation}
when $Q$ is not occupied. With this notation, we define
%
%
%e7.4 #&#
\begin{equation}
\label{energy2} \overline{\mathcal T}_m:={\mathcal T}(X,m)\quad\mbox
{and}\quad\overline\Omega_m:=\bigcup_{Q\in{\mathcal L}(\overline
{\mathcal T}_m)}
\{Q\dvtx \bar\eta_Q >0\}
\end{equation}
to be used in the model selection discussed earlier.

We now describe how to implement the maximization that gives
$\overline{\mathcal T}_m$
and therefore $\overline\Omega_m$. Notice that $\bar\eta
_Q=\gamma_{Q,m}=0$
and ${\mathcal T}(Q,m)$ is empty when $Q$ is not occupied,
and therefore these values are available to us for free.
Thus the computational work in this implementation is solely
determined by the occupied cubes that form ${\mathcal T}'(\mathbf{z})$.
For $l=0,\ldots,\bar n$, we define
%
%
%e7.5 #&#
\begin{equation}
{\mathcal U}_l:={\mathcal T}'(\mathbf{z})\cap{
\mathcal D}_{\bar n-l},
\end{equation}
the set of occupied cubes of resolution level $\bar n-l$. Notice
that ${\mathcal U}_0={\mathcal L}({\mathcal T}'(\mathbf{z}))$.
We work from the leaves of ${\mathcal T}'(\mathbf{z})$ toward the
root, in a
manner similar to CART optimal pruning (see \cite{CART}),
according to the following steps:
\begin{itemize}
\item
$l=0$: we compute for each $Q\in{\mathcal U}_0$ the quantities
$\bar\eta_Q$ and define $\gamma_{Q,0}:= \max\{0,\bar\eta_Q\}$,
${\mathcal T}(Q,0):= \{Q\}$.
This requires at most $\bar n$ arithmetic operations.
\item
For $l=1,\ldots,\bar n$: suppose we have already determined the
quantities $\gamma_{Q,j}$ and $\bar\eta_Q$, as well as the
trees ${\mathcal T}(Q,j)$, for all $Q\in{\mathcal U}_{l-1}$ and $0\le
j\le l-1$. Recall that ${\mathcal T}(Q,j)$ is
a complete subtree. Now for all $0\le j\le l$ and all cubes $Q\in
{\mathcal U}_l$, we compute
%
%
%e7.6 #&#
\begin{equation}
\label{kR} \bigl(\ell^*_j(R)\bigr)_{R\in{\mathcal C}'(Q)}:=
\operatorname{argmax} \biggl\{ \sum_{R\in{\mathcal C}'(Q)}\gamma_{R,\ell
'(R)}\dvtx
\sum_{R\in
{\mathcal C}'(Q)}\ell'(R) =j \biggr\},
\end{equation}
where ${\mathcal C}'(Q):= {\mathcal C}(Q)\cap{\mathcal T}'(\mathbf{z})$
denotes the set of occupied children of $Q$. Notice
that the above argmax may not be unique, in which case we can pick any
maximizer.
We obviously have for each $Q\in{\mathcal U}_l$ and any $1\le j\le l$,
%
%
%e7.7 #&#
\begin{equation}
\label{energy} \gamma_{Q,j}= \sum_{R\in{\mathcal C}'(Q)}
\gamma_{R,\ell^*_{j-1}(R)},
\end{equation}
with
\[
{\mathcal T}(Q,j)= \{Q\}\cup\biggl( \bigcup_{R\in{\mathcal
C}'(Q)}{
\mathcal T}\bigl(R,\ell^*_{j-1}(R)\bigr) \biggr) \cup\bigl({\mathcal
C}(Q)\setminus{\mathcal C}'(Q)\bigr).
\]
For $j=0$, we compute the $\bar\eta_Q$ for all $Q\in{\mathcal
U}_{l}$ by
summing the $\bar\eta_R$ for $R\in{\mathcal C}'(Q)$ and define
$ \gamma_{Q,0}=\max\{0,\bar\eta_Q\}$
and ${\mathcal T}(Q,0)= \{Q\}$.
\item
At the final step $l=\bar n$, the set ${\mathcal U}_{\bar n}$ consists only
of the root $X$, and we have computed ${\mathcal T}(X,m)$ for
$m=0,\ldots,\bar n$.
This provides the estimators $\overline\Omega_m$ for $m=0,\ldots,\bar n$.
\end{itemize}

To estimate the complexity of the algorithm, we need to bound for each
$l \in\{1,\ldots, \bar n\}$
the number of computations required by (\ref{kR}) and (\ref{energy}).
With proper organization, the argmax in (\ref{kR}) can be found
using at most\break ${\mathcal O}(\#(C'(Q))l^2)$ operations.
We can execute (\ref{energy}) with the same order of computation.
The total complexity over all levels is therefore at most ${\mathcal O}(n^4)$
[a~finer analysis can reduce it to~${\mathcal O}(n^3)$].
Also each optimal tree ${\mathcal T}(Q,m)$ can be recorded with at most
$dm$ bits.
It should be noted that the complexity with respect to the data size
$n$ is independent
of the spatial dimension $d$ which only enters when encoding the
optimal trees ${\mathcal T}(X,m)$.

% \begin{remark}
% \label{remd}
%The above estimate exhibits an exponent $2^d$,
%and inspection of the constant $C_d$ shows that it
%also behaves exponentially in $d$. Therefore, the
%complexity of the algorithm in its present form
%is prohibitive in high dimension. This stems
%from the fact that a single refinement gives rise
%to $2^d$ subcubes that are kept in
%the tree even when they are not occupied.
%A much more favorable estimate in the
%high dimensional context is therefore be
%obtained by observing that the sums
%to be computed in \eref{kR} only involve
%the occupied children.
%Another option is to consider
%refinement procedures that correspond to
%binary trees (each cell is split into two).
%These procedures have been described
%for instances for simplicial cells,
% \end{remark}

% \begin{remark}
% \label{remd}

% \end{remark}

%%%%%%%%%%%%%%%%%%%% green end
%%%%%%%%%%%%%%%%%%%%PB

We turn now to the implementation of Algorithm II. We denote by
${\mathcal H}$ the set
of all ($d-1$)-dimensional hyperplanes. Using the notation therein, for
any subtree
${\mathcal T}$ of ${\mathcal T}(\mathbf{z})$ and any $Q\in{\mathcal
L}({\mathcal T})$, the energy is now defined
as
%
%
%e7.8 #&#
\begin{equation}
\label{orn-energy} \gamma_{\mathcal T}:= \sum_{Q\in{\mathcal
L}({\mathcal T})}
\max_{H\in{\mathcal H},i=0,1}\max\{0,\bar\eta_{Q_i^H}\}.
\end{equation}
The set estimator $\overline\Omega_m$ corresponds to a tree
$\overline{\mathcal T}_m\in\frak T_m$ which maximizes the above
energy. Similarly to the previous discussion, we define
%
%
%e7.9 #&#
\begin{equation}
\label{ornate1} \gamma_{Q,0}:= \max_{H\in{\mathcal H},i=0,1}\max\{0,
\bar\eta_{Q_i^H}\}
\end{equation}
and define as
before $\gamma_{Q,m}$ and ${\mathcal T}(Q,m)$ by (\ref{energy1}) and
(\ref{energy2}).

The procedure of determining the trees ${\mathcal T}(X,m)$ for
$m=0,\ldots,k$ is then, in principle, the same as above, however with
a significant distinction due to the
search for a ``best'' hyperplane
$H=H_Q$ that attains the maximum in (\ref{ornate1}).
Since a cube $Q$ contains a finite number
$n_Q$ of data, the search can be
reduced to ${n_Q\choose d}$ hyperplanes
and the cost of computing $\gamma_{Q,0}$
is therefore bounded by $n_Q^{d}$.
In addition, the search of $H_Q$ needs
to be performed on \emph{every} cube $Q\in{\mathcal T}(\mathbf{z})$,
so that a crude global bound for this cost is given
by $n^{d+2}$. This additional cost is affordable
for small $d$ but becomes prohibitive in
high dimension. An alternate strategy is
to rely on more affordable classifiers to produce an affine (or even
higher order algebraic)
decision boundary on each $Q$.
Examples are plug-in classifiers that
are based on least-square estimation of $\eta$ on $Q$ by a polynomial.

\section*{Acknowledgments}
The authors wish to thank Stephane Gaiffas {{and L\'{a}szl\'{o}
Gy\"orfi}} for various valuable suggestions and references, as well as the
anonymous referees for their constructive comments.

\begin{supplement}[id=suppA]
\stitle{Proof of Theorem 2.1}
\slink[doi]{10.1214/14-AOS1234SUPP} %[doi,text={...}] - jei reikia suskaldyti doi
\sdatatype{.pdf}
\sfilename{aos1234\_supp.pdf}
\sdescription{This supplement contains the detailed proof of Theorem~2.1 \cite{BCDD3}.}
\end{supplement}

% imsref loaded by linak, 2014-07-25 16:25:45
%
% imsref loaded by linak, 2014-07-30 12:00:04
% imsref loaded by linak, 2014-09-18 09:47:56

\printaddresses
\end{document}